\documentclass[reqno,12pt]{amsart}
\usepackage{amsmath,amsthm,amssymb,amsfonts,amscd}
\usepackage{xypic}
\theoremstyle{plain} 
	\newtheorem{thm}{Theorem}[section]
	\newtheorem*{thm*}{Theorem}
	\newtheorem{cor}[thm]{Corollary}
	\newtheorem{lem}[thm]{Lemma}
	
	\newtheorem{prop}[thm]{Proposition}
	
	\newtheorem*{conj*}{Conjecture}
	
\theoremstyle{definition}
	\newtheorem{defn}[thm]{\rm{Definition}}

\theoremstyle{remark}
	\newtheorem{rmk}[thm]{\rm{Remark}}
	
	\newtheorem*{pf}{\rm{Proof}}
\numberwithin{equation}{section}
\def\CC{{\mathbb C}}

\def\RR{{\mathbb R}}
\def\ZZ{{\mathbb Z}}

\def\p{\partial }

\def\multi{{\mathbb{Z}^n_{\geq 0}}}
\def\one{{1 \leq \alpha \leq n}}
\def\two{{1 \leq \alpha,\,\beta \leq n}}
\def\Der{{Der(\CC[V]^G)}}

\newcommand{\bp}{\begin{pmatrix}}
\newcommand{\ep}{\end{pmatrix}}

\numberwithin{equation}{section}

\newcounter{CounterEQUlabel}
\newcommand{\EQUlabel}[1]{\label{#1}
	\ifcase \theCounterEQUlabel
		\relax
	\or
		\hspace{1em}\mbox{\tiny$\langle$\rmfamily#1$\rangle$}
		\index{zzz#1@#1}
	\fi }	
	\newcounter{CounterEQUref}
	\newcounter{CounterEQUpageref}
	\newcommand{\EQUref}[1]{
		\ifcase \theCounterEQUref     \relax   \or {\tiny[#1]}\,\fi
		\ifcase \theCounterEQUpageref (\ref{#1}) \or (\ref{#1}\,(p.\pageref{#1})) \fi}

\begin{document}
\title{Good basic invariants and Frobenius structures}
\date{\today}
\author{Ikuo Satake}
\address{Faculty of Education, Bunkyo University, 
3337 Minamiogishima Koshigaya, Saitama, 343-8511, Japan}
\email{satakeikuo@gmail.com}
\begin{abstract}
In this paper, we define a set of good basic invariants for 
a finite complex reflection group under certain conditions. 
We show that a set of good basic invariants 
for a finite real reflection group 
gives a set of the flat invariants obtained by Saito and 
the Taylor coefficients of these good basic invariants 
give the structure constants of the multiplication of the 
Frobenius structure obtained by Dubrovin. 
\end{abstract}
\maketitle
\section{Introduction}

\subsection{Aim and results of the paper}
Let $G$ be a finite real reflection group which acts on the real vector space 
$V_{\RR}$, 
$g \in G$ be a Coxeter transformation and 
$q \in V_{\RR}\otimes_{\RR}\CC$ be an eigenvector 
of $g$ whose eigenvalue $\zeta$ is 
a primitive $h$-th root of unity, where 
$V_{\RR}\otimes_{\RR}\CC$ is the complexification of the vector space $V_{\RR}$ and $h$ is the Coxeter number. 
It is known that the eigenvector $q$ is regular 
(that is, it does not lie on any reflecting hyperplane). 

In the theory of the flat structure (the Frobenius structure)
(see \cite{finite}) for the invariants of the finite real reflection group, 
the values of $G$-invariants and $G$-anti-invariants 
at the regular eigenvector $q$ play an important role. 

In this paper, we first assume that 
$G$ is a finite complex reflection group. 
Then we study the Taylor expansions 
(with a suitable grading) of $G$-invariants 
at a suitable regular vector $q$ 
and define a set of ``good basic invariants" 
by using Taylor expansions at $q$ under certain conditions 
(the existence of an admissible triplet defined in Definition \ref{200310.1}). 
If $G$ is a finite real reflection group, we show that 
\begin{enumerate}
\item a set of good basic invariants gives a set of
 flat invariants obtained by Saito \cite{finite} of the Frobenius structure,
\item the Taylor coefficients of the good basic invariants 
give the structure constants of the multiplication of the Frobenius structure 
obtained by Dubrovin \cite{Dubrovin}. 
\end{enumerate}

Here is a brief account of the contents of the paper. 
In Section 2 we define an admissible triplet for 
the finite complex reflection group. 
In Section 3 we define good basic invariants. 
In Section 4 we give examples of good basic invariants. 
In Section 5 we study the dependence of good basic invariants 
on the choice of the admissible triplet. 
In Section 6 we give properties of Taylor coefficients of 
good basic invariants. 
In Section 7 we treat the cases of finite real reflection groups. 
We show the existence and the uniqueness of good basic invariants 
and we give a description of the bilinear form in terms of 
the good basic invariants (Theorem \ref{230403.3}). 
In Section 8 we show that the good invariants give a nice description 
of the Frobenius structure which is defined by 
Saito and Dubrovin. 

\subsection{Acknowledgements}
The author thanks Prof. Yukiko Konishi and Prof. 
Satoshi Minabe for the careful reading of the manuscript. 

This work is supported in part by 
Grant-in Aid for Challenging Research (Exploratory) 17K18781
, 
Grant-in-Aid for Scientific Research(C) 18K03281 
and 
Grant-in-Aid for Scientific Research(C) 22K03295.  
\\

\section{Graded $\CC$-algebra structure on $\mathbb{C}[V]$}
Let $V$ be a $\mathbb{C}$-vector space of $\mathrm{dim}_{\mathbb{C}}V=n$. 
From Section 2 to Section 5, we assume that 
$G \subset GL(V)$ is a finite complex reflection group, i.e. 
$G$ is generated by reflections, where $g \in GL(V)$ is called 
a reflection if it is of finite order and if 
all but one of its eigenvalues are equal to $1$. 
We also assume that $G$ is irreducible, i.e. $V$ is an irreducible $G$-module.

\subsection{Graded $\CC$-algebra structure on $\mathbb{C}[V]^G$}
Let $\mathbb{C}[V]$ be a symmetric tensor algebra 
of $V^*:=\mathrm{Hom}(V,\mathbb{C})$, which is identified with 
the algebra of polynomial functions on $V$. 
Let $z^1,\cdots,z^n$ be a basis of $V^*$. 
Then a set of 
\begin{equation}
z^a:=(z^1)^{a_1}\cdots (z^n)^{a_n}\quad
(a=(a_1,\cdots,a_n) \in \mathbb{Z}^n_{\geq 0})
\end{equation}
gives a $\mathbb{C}$-basis of $\mathbb{C}[V]$. 

On the $\mathbb{C}$-algebra $\mathbb{C}[V]$, the natural grading 
is defined by counting the degree of $z^a$ as 
\begin{equation}
|a|=a_1+\cdots+a_n 
\end{equation}
for $a=(a_1,\cdots,a_n) \in \mathbb{Z}^n_{\geq 0}$. 

The action of $g \in G$ on $F \in \mathbb{C}[V]$ is defined by 
\begin{equation}
(g\cdot F)(v)=F(g^{-1}\cdot v)\quad
(v \in V). 
\end{equation}
We denote the algebra of $G$-invariant elements of $\mathbb{C}[V]$ 
by $\mathbb{C}[V]^G$. 
The grading of $\mathbb{C}[V]$ induces the grading on 
the algebra $\mathbb{C}[V]^G$. 
We denote its degree $j$ part by $S(j)$. 

By the famous result of Shephard-Todd-Chevalley, 
for a finite complex reflection group $G$,  
the algebra $\mathbb{C}[V]^G$ is 
generated by algebraically independent 
homogeneous elements 
$x^1,\cdots,x^n$ with degree $d_1 \leq d_2 \leq \cdots \leq d_n$ 
which we call a set of {\it basic invariants}. 

We remark that a set of basic invariants is not unique, 
but the degrees $d_1,\cdots,d_n$ are uniquely determined. 
We put 
\begin{equation}
d:=(d_1,\cdots,d_n).
\end{equation}

Then the degree $j$ part $S(j)$ could be also written as 
\begin{equation}
S(j):=\{\sum_{b \in \multi} A_b x^b\in \mathbb{C}[V]^G\,|\, A_b \in \mathbb{C},\ 
d\cdot b=j\}
\end{equation}
where we denote 
\begin{equation}
x^b=(x^1)^{b_1}\cdots (x^n)^{b_n},\quad
d\cdot b=d_1 b_1+\cdots+d_nb_n 
\end{equation}
and we have the decomposition 
\begin{equation}\label{200309.1}
\mathbb{C}[V]^G=\bigoplus_{j=0}^{\infty}S(j).
\end{equation}
\subsection{Admissible triplet}
In this subsection, we introduce the notion of an admissible triplet. 
\begin{defn}\label{200310.1}
For $g \in G,\ \zeta \in \mathbb{C}$ and $q \in V$, 
we call a triplet $(g,\zeta,q)$ {\it admissible} if it 
satisfies the following conditions. 
\begin{enumerate}
\item the vector $q \in V$ is an eigenvector of $g$ with the eigenvalue 
$\zeta$ which is a primitive $d_n$-th root of unity.
\item the Jacobian matrix 
\begin{equation}\label{200318.10}
\left(
\frac{\p x^{\alpha}}{\p z^{\beta}}(q)
\right)_{1 \leq \alpha,{\beta} \leq n}
\end{equation}
is invertible, where $z^1,\cdots,z^n \in V^*$ form a basis of $V^*$. 
\end{enumerate}
\end{defn}

In this section, we fix an admissible triplet $(g,\zeta,q)$. 
By the action of $g$ on the Jacobian matrix (\ref{200318.10}), 
we see that the eigenvalues of $g$ on $V$ are $\zeta^{1-d_{\alpha}}$ 
$(1 \leq \alpha \leq n)$ (cf. Theorem 4.2\hskip0.2mm(v) of [5]). 
Hence we may and shall assume that 
\begin{equation}\label{200316.6}
g \cdot z^{\alpha}=\zeta^{d_{\alpha}-1}z^{\alpha}\ (1 \leq {\alpha} \leq n)
\end{equation}
for $z^{\alpha}$ in Definition \ref{200310.1}\hskip0.2mm(ii). 
A basis of $V^*$ satisfying (\ref{200316.6}) is called a ``$g$-homogeneous basis". 
We put 
\begin{equation}
\widetilde{g}:=\zeta^{-1}\cdot g \in GL(V).
\end{equation}
Then we have 
\begin{equation}\label{200316.2}
\widetilde{g} \cdot q=q,\quad
\widetilde{g} \cdot z^{\alpha}=\zeta^{d_{\alpha}}z^{\alpha},\quad
\widetilde{g} \cdot x^{\alpha}=\zeta^{d_{\alpha}}x^{\alpha}\ (1 \leq {\alpha} \leq n).
\end{equation}
Thus we have 
\begin{eqnarray}
&&z^{\alpha}(q)=0\ (d_{\alpha}<d_n),\quad
x^{\alpha}(q)=0\ (d_{\alpha}<d_n), \label{200310.2}\\
&&\frac{\p x^{\alpha}}{\p z^{\beta}}(q)=0\ 
(d_{\alpha}\neq d_{\beta}).
\label{200316.1}
\end{eqnarray}
For $a,b \in \multi$, we have 
\begin{equation}\label{230402.1}
\left(
\frac{\p^b}{\p z^b}x^a
\right)
(q)=0\quad 
(d\cdot b \not\equiv d\cdot a\ (\mathrm{mod}\ d_n)), 
\end{equation}
where we denote 
\begin{equation}
\frac{\p^b}{\p z^b}=
\left(\frac{\p}{\p z^1}\right)^{b_1}
\cdots
\left(\frac{\p}{\p z^n}\right)^{b_n}\quad
\mbox{ for }b=(b_1,\cdots,b_n) \in \multi. 
\end{equation}

\begin{rmk}
We remark that the facts (\ref{200318.10}), (\ref{200316.2}), (\ref{200310.2}), 
(\ref{200316.1}) suggest that 
we should see $q$ as a fixed point of $\widetilde{g}$ and also 
see a $G$-invariant element $x^{\alpha}$ as an analogous object of $z^{\alpha}$ at $q$ 
for $1 \leq \alpha \leq n$. 
This leads to the idea to study a Taylor expansion of 
$x^{\alpha}-x^{\alpha}(q)$ $(1 \leq \alpha \leq n)$ 
by $z^1-z^1(q),\cdots,z^n-z^n(q)$. 
\end{rmk}
%
%
%
\subsection{Graded $\CC$-algebra structure on $\mathbb{C}[V]$}
We introduce on $\mathbb{C}[V]$ 
another $\mathbb{Z}$-grading by the aid of the admissible triplet $(g,\zeta,q)$. 
\begin{defn}
Let $z^1,\cdots,z^n$ be a $g$-homogeneous basis of $V^*$. 
For any $j \in \mathbb{Z}_{\geq 0}$, we define 
\begin{equation}\label{200316.7}
V(g,\zeta,q)(j):=
\{\sum_{b \in \multi} c_b z^b\in \mathbb{C}[V]\,|\,c_b \in \CC, d\cdot b=j\}.
\end{equation}
\end{defn}
We remark that the admissible triplet $(g,\zeta,q)$ naturally gives 
only $\mathbb{Z}/d_n\mathbb{Z}$-grading 
on $\mathbb{C}[V]$ by $\widetilde{g} \cdot z^{\alpha}=\zeta^{d_{\alpha}}z^{\alpha}$ 
for $\one$. 
We lift it to the $\mathbb{Z}$-grading.

We give a graded $\CC$-algebra structure on $\mathbb{C}[V]$ 
by the decomposition 
\begin{equation}\label{200309.2}
\mathbb{C}[V]=\bigoplus_{j=0}^{\infty}V(g,\zeta,q)(j).
\end{equation}
%
%
\\

\section{Graded $\CC$-algebra Isomorphism $\psi$}
Let $G$ be a finite complex reflection group and 
we fix an admissible triplet $(g,\zeta,q)$. 
\subsection{The morphism $\varphi[g,\zeta,q]$}
\begin{defn}\label{230403.1}
For the admissible triplet $(g,\zeta,q)$, 
a set of basic invariants $x^1,\cdots,x^n$ 
and a $g$-homogeneous basis $z^1,\cdots,z^n$ of $V^*$, 
we define a $\mathbb{C}$-module homomorphism:
\begin{equation}
\varphi[g,\zeta,q]:\mathbb{C}[V]^G \to \mathbb{C}[V],\quad
x^a \mapsto \sum_{b \in \multi}
\frac{1}{b!}
\frac{\partial^b([x-x(q)]^a)}{\partial z^b}(q)z^b,\label{200310.13}
\end{equation}
for a $\mathbb{C}$-basis $\{x^a\,|\,a \in \multi\}$ of 
$\mathbb{C}$-module $\mathbb{C}[V]^G$, where 
we used notations
\begin{eqnarray}
[x-x(q)]^a&:=&(x^1-x^1(q))^{a_1}\cdots (x^n-x^n(q))^{a_n} 
\ \mbox{for}\ a=(a_1,\cdots,a_n) \in \ZZ^n_{\geq 0},\nonumber \\ 
b!&:=&b_1!\cdots b_n!\ \mbox{for}\ b=(b_1,\cdots,b_n) \in \ZZ^n_{\geq 0}.
\end{eqnarray}
\end{defn}
We remark that $\varphi[g,\zeta,q](f)$ for $f \in \CC[V]^G$ 
is not necessarily invariant by the $G$-action. 

\begin{prop}\label{200310.12}
\begin{enumerate}
\item $\varphi[g,\zeta,q]$ depends neither on the choices of 
a set of basic invariants $x^{\alpha}$ 
$(\one)$ nor on the choices of a $g$-homogeneous basis $z^{\alpha}$ $(\one)$ of $V^*$. 
$\varphi[g,\zeta,q]$ gives a $\mathbb{C}$-algebra homomorphism. 
\item Let $x^{\alpha}\,(1 \leq \alpha \leq n)$ 
and $z^{\alpha}\,(1 \leq \alpha \leq n)$ be the same as 
in Definition $\ref{230403.1}$. 
For any multi-indices $a,b \in \multi$, the coefficients of $z^b$ 
of the RHS of 
$$
\varphi[g,\zeta,q](x^a)
=
\sum_{b \in \multi}
\frac{1}{b!}
\frac{\partial^b([x-x(q)]^a)}{\partial z^b}(q)z^b
$$
is $0$ if $d\cdot b \notin \{d\cdot a+d_n j\,|\,j \in \ZZ_{\geq 0}\}$. 
\end{enumerate}
\end{prop}
\begin{pf} 
(i) For a set of basic invariants $x^{\alpha}\ (1 \leq \alpha \leq n)$, 
we define a $\mathbb{C}$-algebra homomorphism $\varphi_1[q]$ by 
\begin{equation}
\varphi_1[q]:\mathbb{C}[V]^G \to \mathbb{C}[V]^G,\quad 
x^{a} \mapsto (x-x(q))^a\ (a \in \multi). 
\end{equation}
Let $k$ be an integer satisfying $d_{k}<d_n$ and $d_{k+1}=d_n$. 
Put $S_0=\mathbb{C}[x^1,\cdots,x^{k}]$. 
Then $\mathbb{C}[V]^G$ is an $S_0$-algebra with polynomial generators 
$x^{k+1},\cdots,x^n$. 
The $\mathbb{C}$-algebra homomorphism $\varphi_1[q]$ leaves invariant 
an element of $S_0$ 
because the basic invariants $x^{\alpha}\ (d_{\alpha} <d_n)$ 
satisfies $x^{\alpha}(q)=0$ by (\ref{200310.2}). 
Then $\varphi_1[q]$ is an $S_0$-algebra homomorphism which is determined 
by $\varphi_1[q](x^{k+1}),\cdots,\varphi_1[q](x^n)$.


If we take another set of basic invariants $y^1,\cdots,y^n$, 
then $\mathbb{C}[y^1,\cdots,y^k]=S_0$ and 
$y^{k+1},\cdots,y^{n}$ are sums of $\mathbb{C}$-linear combinations of 
$x^{k+1},\cdots,x^n$ and elements of $S_0$, i.e. 
for $k+1 \leq j \leq n$, we have 
$$
y^j=\sum_{i=k+1}^n c_{i}^{j}x^i
+\sum_
{\substack{
a=(a_1,\cdots,a_n), d\cdot a=d_j, \\
a_{k+1}=\cdots=a_n=0}}
c_{a}^j x^a, 
$$
where $c_{i}^{j},\ c_{a}^j \in \CC$. 
Then we have 
\begin{equation*}
y^j-y^j(q)=
\sum_{i=k+1}^n c_{i}^{j}(x^i-x^i(q))+
\sum_
{\substack{
a=(a_1,\cdots,a_n), d\cdot a=d_j, \\
a_{k+1}=\cdots=a_n=0}}
c_{a}^j (x^a-x^a(q)).
\end{equation*}
For any $(x^a-x^a(q))$ in the second sum, we have 
$(x-x(q))^a=(x^a-x^a(q))$ by $x^{\gamma}(q)=0$ for all $\gamma \leq k$. 
Thus we have 
\begin{equation}\label{240307.1}
y^j-y^j(q)=
\sum_{i=k+1}^n c_{i}^{j}(x^i-x^i(q))+
\sum_
{\substack{
a=(a_1,\cdots,a_n), d\cdot a=d_j, \\
a_{k+1}=\cdots=a_n=0}}
c_{a}^j (x-x(q))^a.
\end{equation}
This means that the morphism $\varphi_1[q]$ does not depend on the choice of a set of basic invariants 
$x^1,\cdots,x^n$ but depends only on the choice of $q$.

We define a $\mathbb{C}$-algebra homomorphism $\varphi_2$ by 
\begin{equation}
\varphi_2:\mathbb{C}[V]^G \to \mathbb{C}[V],\quad
f \mapsto \sum_{b \in \multi}\frac{1}{b!}
\frac{\partial^b f}{\partial z^b}(q)(z-z(q))^b. 
\end{equation}
This is a Taylor expansion at $q$ and it coincides with the natural inclusion 
$\mathbb{C}[V]^G \subset \mathbb{C}[V]$. 

We define a $\mathbb{C}$-algebra homomorphism $\varphi_3[q]$ by 
\begin{equation}
\varphi_3[q]:\mathbb{C}[V] \to \mathbb{C}[V], 
\quad
z^{a} \mapsto (z+z(q))^a\quad
(a \in \multi). 
\end{equation}
This morphism does not depend on the choice of a $g$-homogeneous basis 
$z^1,\cdots,z^n$ of $V^*$ because a basis is unique up to linear transformations. 

Then we have 
\begin{equation}
\varphi[g,\zeta,q]=\varphi_3[q] \circ \varphi_2 \circ \varphi_1[q]. 
\end{equation}
Since $\varphi_1[q]$, $\varphi_2$ and $\varphi_3[q]$ are $\mathbb{C}$-algebra homomorphisms, 
their composite morphism $\varphi[g,\zeta,q]$ is also 
a $\mathbb{C}$-algebra homomorphism. 

(ii) If the assersion (ii) is true for the multi-indices $a,a' \in \multi$, 
then it is true for the multi-index $a+a' \in \multi$. 
Thus we should only prove (ii) for each element of a set 
of the basic invariants $x^1,\cdots,x^n$. 

For any $x^{\alpha}\ (\one)$, we have
\begin{eqnarray*}
\varphi[g,\zeta,q](x^{\alpha})&=&
\sum_{b \in \multi}\frac{1}{b!}
\frac{\partial^b(x^{\alpha}-x^{\alpha}(q))}{\partial z^b}(q)z^b\\
&=&
(x^{\alpha}-x^{\alpha}(q))(q)+
\sum_{b \in \multi,\ b\neq 0}\frac{1}{b!}
\frac{\partial^b(x^{\alpha}-x^{\alpha}(q))}{\partial z^b}(q)z^b\\
&=&\sum_{b \in \multi,\ b\neq 0}
\frac{1}{b!}
\frac{\partial^b x^{\alpha}}{\partial z^b}(q)z^b. 
\end{eqnarray*}
The coefficients of $z^b$ are $0$ if 
$d \cdot b \notin \{d_{\alpha}+d_n j\,|\,j \in \ZZ_{\geq 0}\}$ 
by (\ref{230402.1}). This gives a proof of (ii). 
\qed\end{pf}
\subsection{The morphism $\psi[g,\zeta,q]$}

In this subsection, we use properties of a filtered algebra 
(cf. \cite[Ch.3 \S 2-3]{commutative}). 
The $\CC$-algebra homomorphism $\varphi[g,\zeta,q]$ 
is not a graded $\CC$-algebra homomorphism 
with respect to the grading (\ref{200309.1}) on $\CC[V]^G$ 
and (\ref{200309.2}) on $\CC[V]$. 

However if we define decreasing filtrations on $\CC[V]^G$ 
and $\CC[V]$ by 
\begin{eqnarray}
F^{m}(\CC[V]^G)&:=&
\bigoplus_{j \geq m}S(j)
\quad (\forall m \in \ZZ_{\geq 0}),
\\
F^{m}(\CC[V])&:=&
\bigoplus_{j \geq m}V(g,\zeta,q)(j)
\quad
(\forall m \in \ZZ_{\geq 0})
\end{eqnarray}
respectively, $\CC[V]^G$ and $\CC[V]$ are filtered $\CC$-algebras 
and $\varphi[g,\zeta,q]$ is a filtered $\CC$-algebra homomorphism 
because we have 
\begin{equation}
\varphi[g,\zeta,q](F^{m}(\CC[V]^G))\subset F^{m}(\CC[V])
\quad (\forall m \in \ZZ_{\geq 0})
\end{equation}
by Proposition \ref{200310.12} (ii). 
\begin{defn}\label{230804.1}
\begin{enumerate}
\item 
Let $\mathrm{gr}_F\varphi[g,\zeta,q]$ be the 
graded $\CC$-algebra homomorphism induced by 
a filtered $\CC$-algebra homomorphism $\varphi[g,\zeta,q]$: 
\begin{equation}
\mathrm{gr}_F\varphi[g,\zeta,q]:\mathrm{gr}_F(\CC[V]^G) 
\to \mathrm{gr}_F(\CC[V]),
\end{equation}
where 
\begin{eqnarray}
\mathrm{gr}_F(\CC[V]^G)
&:=&\bigoplus_{m \in \ZZ_{\geq 0}}F^{m}(\CC[V]^G)/F^{m+1}(\CC[V]^G),\\
\mathrm{gr}_F(\CC[V])
&:=&\bigoplus_{m \in \ZZ_{\geq 0}}F^{m}(\CC[V])/F^{m+1}(\CC[V]). 
\end{eqnarray}
\item Let $\psi[g,\zeta,q]$ be the graded $\CC$-algebra homomorphism defined by
\begin{equation}
\psi[g,\zeta,q]:=\psi_2^{-1} \circ \mathrm{gr}_F\varphi[g,\zeta,q] 
\circ \psi_1:
\CC[V]^G \to \CC[V], 
\end{equation}
where we used the natural graded $\CC$-algebra isomorphism
\begin{eqnarray}
&&\psi_1:\CC[V]^G \to \mathrm{gr}_F(\CC[V]^G) \\
&&(resp.\,\psi_2:\CC[V] \to \mathrm{gr}_F(\CC[V])), 
\end{eqnarray}
which maps an element of $S(j)$ (resp.\hskip1mm$V(g,\zeta,q)$) 
to its canonical image in \\
$F^j(\CC[V]^G)/F^{j+1}(\CC[V]^G)$ (resp.\hskip1mm$F^{j}(\CC[V])/F^{j+1}(\CC[V])$). 
\end{enumerate}
\end{defn}
We have an explicit description of $\psi[g,\zeta,q]$:
\begin{equation}
\psi[g,\zeta,q]:\mathbb{C}[V]^G \to \mathbb{C}[V],\quad
x^a \mapsto \sum_{b \in \multi,\ 
d\cdot b=d\cdot a}\frac{1}{b!}
\frac{\partial^b([x-x(q)]^a)}{\partial z^b}(q)z^b \label{200310.14}
\end{equation}
for a $\mathbb{C}$-basis $\{x^a\,|\,a \in \multi\}$ of 
$\mathbb{C}$-module $\mathbb{C}[V]^G$, 
where we used notations in Definition \ref{230403.1}.

\begin{prop}
With respect to the gradings $($\ref{200309.1}$)$ on $\mathbb{C}[V]^G$ 
and $($\ref{200309.2}$)$ on $\mathbb{C}[V]$, 
$\psi[g,\zeta,q]$ is a graded $\CC$-algebra isomorphism
\begin{equation}
\psi[g,\zeta,q]:
\CC[V]^G
\stackrel{\sim}{\to}
\CC[V].
\end{equation}
\end{prop}
\begin{pf}
In our proof, we denote $\psi[g,\zeta,q]$ simply by $\psi$. 

For a proof, we have only to prove that 
$\{\psi(x^{1}),\cdots,\psi(x^{n})\}$ is a set of homogeneous polynomial 
generators of the graded $\CC$-algebra $\mathbb{C}[V]$. 
%
It is equivalent to show that 
$\psi(x^{\alpha})\ (1 \leq \alpha \leq n)$ 
is written as 
\begin{equation}\label{210721.1}
\psi(x^{\alpha})=\sum_{1 \leq \beta \leq n,\, d_{\beta}=d_{\alpha}}
A^{\alpha}_{\beta}z^{\beta}
+\sum_{a \in \multi,\, d\cdot a=d_{\alpha},|a|\geq 2}B^{\alpha}_a z^a
\end{equation} 
with the matrix $(A^{\alpha}_{\beta})$ invertible, where we put 
$A^{\alpha}_{\beta}=0$ for 
$d_{\alpha} \neq d_{\beta}\ (1 \leq \alpha,\beta \leq n)$.

By the admissibility (ii) of the triplet $(g,\zeta,q)$, 
the Jacobian matrix 
\begin{equation}
\left(
\frac{\p x^{\alpha}}{\p z^{\beta}}(q)
\right)_{1 \leq \alpha,\beta \leq n}
\end{equation}
is invertible. 
For any $\alpha,\beta\ (1 \leq \alpha,\beta \leq n)$, the entry 
$$
\frac{\p x^{\alpha}}{\p z^{\beta}}(q)
$$ 
of the Jacobian matrix is $0$ 
if $d_{\alpha} \neq d_{\beta}$ by (\ref{200316.1}). 

Then the Jacobian matrix $J$ is a block diagonal matrix with each block  
\begin{equation}\label{200316.8}
\left(
\frac{\p x^{\alpha}}{\p z^{\beta}}(q)
\right)_{d_{\alpha}=d_{\beta}=k}
\end{equation}
which is an invertible matrix. 

Then 
\begin{equation}
\psi(x^{\alpha})=\sum_{b \in \multi,\,d\cdot b=d_{\alpha}}\frac{1}{b!}
\frac{\partial^b[x^{\alpha}-x^{\alpha}(q)]}{\partial z^b}(q)z^b \ (1 \leq \alpha \leq n)
\end{equation}
satisfies the conditions (\ref{210721.1}) and 
we see that 
$\{\psi(x^{1}),\cdots,\psi(x^{n})\}$ 
gives a set of homogeneous polynomial 
generators of the graded $\CC$-algebra $\mathbb{C}[V]$. 
\qed\end{pf}

\subsection{Good basic invariants}
\begin{defn}\label{200316.9}
A set of basic invariants $x^1,\cdots,x^n$ is good with respect to 
the admissible triplet $(g,\zeta,q)$ if 
$x^1,\cdots,x^n$ form a $\mathbb{C}$-basis of the vector space 
$\psi[g,\zeta,q]^{-1}(V^*)$ w.r.t. the natural inclusion $V^* \subset \mathbb{C}[V]$. We call $x^1,\cdots,x^n$ ``good basic invariants". \\
\end{defn}



\section{Examples}
In this section, we give some examples of a set of 
good basic invariants. 
 
Let $\CC^{l+1}$ be a $\CC$-vector space with coordinates 
$\varepsilon^1,\cdots,\varepsilon^{l+1}$ 
and 
the symmetric group $S_{l+1}$ acts on $\CC^{l+1}$ by 
$\sigma(\varepsilon^1,\cdots,\varepsilon^{l+1}) 
=(\varepsilon^{\sigma^{-1}(1)},\cdots,
\varepsilon^{\sigma^{-1}(l+1)})$ for $\sigma \in S_{l+1}$. 
This action preserves the subspace $V$ 
where 
$V=\{(\varepsilon^1,\cdots,\varepsilon^{l+1})\in \CC^{l+1}\,|\,
\sum_{i=1}^{l+1}\varepsilon^i=0\}$. 
Then this action gives an injection $S_{l+1} \to GL(V)$. 
By this injection, we regard the group $S_{l+1}$ 
as a finite complex reflection group. 

We define $g \in S_{l+1}$ by 
$g(i)=i-1$ for $i=2,\cdots,l+1$ and $g(1)=l+1$. 
We put 
$q=(1,\zeta,\cdots,\zeta^l) \in V$, 
where $\zeta=\exp(\frac{2\pi\sqrt{-1}}{l+1})$. 
Then we have $g\cdot q=\zeta q$ and $(g,\zeta,q)$ gives 
an admissible triplet. 
We remark that this is the case of $A_l$-type 
(cf. Proposition \ref{230403.2}). 

We define linear functions on $\CC^{l+1}$ by 
$$
y^k=\sum_{i=1}^{l+1} \zeta^{(i-1)(k-1)}\varepsilon^i\quad (k=1,\cdots,l+1). 
$$
Then we have $g\cdot y^k=\zeta^{k-1}y^k$ for $k=1,\cdots,l+1$. 

We define $S_{l+1}$-invariant polynomial functions on $\CC^{l+1}$ by  
$$
P^1=(-1)^{1-1}\sum_{i=1}^{l+1} \varepsilon^i, \ 
P^2=(-1)^{2-1}\sum_{1 \leq i<j \leq l+1}^{l+1} 
\varepsilon^i\varepsilon^j, \ 
\cdots,
P^{l+1}=(-1)^{(l+1)-1}\varepsilon^1\cdots\varepsilon^{l+1}. 
$$

Then $y^2,\cdots,y^{l+1}$ give a basis of $V^*$ and 
$P^2,\cdots,P^{l+1}$ give a set of basic invariants of $\CC[V]^{S_{l+1}}$. 

We give an explicit description of $\varphi[g,\zeta,q]$ and 
$\psi[g,\zeta,q]$ for $l=1,\,2,\,3$ cases. 
We remark that the degree of 
$y^a=(y^2)^{a_2}\cdots (y^{l+1})^{a_{l+1}}$ for 
$a=(a_2,\cdots,a_{l+1}) \in \ZZ^{l}_{\geq 0}$ 
is $\sum_{k=2}^{l+1} ka_{k}$, i.e. 
$$
y^a \in V(g,\zeta,q)(\sum_{k=2}^{l+1} ka_{k}). 
$$
We use below square brackets in order to put together 
the same degree terms. 

$l=1$ case. We have 
\begin{eqnarray*}
P^2&=&\frac{1}{4}(y^2)^2,\\
\varphi[g,\zeta,q](P^2)&=&y^2+\frac{1}{4}(y^2)^2, \\
\psi[g,\zeta,q](P^2)&=&y^2. 
\end{eqnarray*}
Then $P^2$ give a set of good basic invarant. 

$l=2$ case. We have 
\begin{eqnarray*}
P^2&=&\frac{1}{3}y^2y^3,\\
P^3&=&\frac{1}{3^3}((y^3)^3+(y^2)^3),\\
\varphi[g,\zeta,q](P^2)&=&y^2+\frac{1}{3}y^2y^3,\\
\varphi[g,\zeta,q](P^3)&=&y^3+\left[\frac{1}{3}(y^3)^2
+\frac{1}{3^3}(y^2)^3\right]
+\frac{1}{3^3}(y^3)^3,\\
\psi[g,\zeta,q](P^2)&=&y^2,\\
\psi[g,\zeta,q](P^3)&=&y^3.
\end{eqnarray*}
Then $P^2,\ P^3$ give a set of good basic invarants. 

$l=3$ case. We have 
\begin{eqnarray*}
P^2&=&\frac{1}{4^2}\left[4 y^2y^4+2(y^3)^2\right],\\
P^3&=&\frac{1}{4^3}\left[4y^3(y^4)^2+4(y^2)^2y^3\right],\\
P^4&=&\frac{1}{4^4}\left[(y^4)^4-2(y^2)^2(y^4)^2+4y^2(y^3)^2y^4-(y^3)^4+(y^2)^4\right],\\
\varphi[g,\zeta,q](P^2)&=&y^2+\left[\frac{1}{4}y^2y^4+\frac{1}{8}(y^3)^2\right],\\
\varphi[g,\zeta,q](P^3)&=&y^3+\left[\frac{1}{4}(y^2)^2y^3+\frac{1}{2}y^3y^4\right]
+\frac{1}{4^2}y^3(y^4)^2,\\
\varphi[g,\zeta,q](P^4)&=&\left[y^4-\frac{1}{8}(y^2)^2\right]
+\left[\frac{6}{4^2}(y^4)^2
-\frac{1}{4^2}(y^2)^2y^4
+\frac{1}{4^2}y^2(y^3)^2
+\frac{1}{4^4}(y^2)^4\right]\\
&&
+\left[\frac{1}{4^2}(y^4)^3
-\frac{2}{4^4}(y^2)^2(y^4)^2
+\frac{1}{4^3}y^2(y^3)^2y^4
-\frac{1}{4^4}(y^3)^4\right]
+\frac{1}{4^4}(y^4)^4,\\
\psi[g,\zeta,q](P^2)&=&y^2,\\
\psi[g,\zeta,q](P^3)&=&y^3,\\
\psi[g,\zeta,q](P^4)&=&y^4-\frac{1}{8}(y^2)^2.
\end{eqnarray*}
Then $P^2,\ P^3,\ P^4+\frac{1}{8}(P^2)^2$ give a set of good basic invarants. 
\\

\section{Independence on the choice of the admissible triplet}
Let $G$ be a finite complex reflection group. 
We study the dependence of good basic invariants 
on the choice of the admissible triplet $(g,\zeta,q)$.  

For any $h \in GL(V)$ and a $g$-homogeneous basis $z^1,\cdots,z^n$ of $V^*$, we define 
\begin{equation}
\xi[h]:\mathbb{C}[V] \to \mathbb{C}[V], \quad
z^b \mapsto (h\cdot z)^b \quad 
(b \in \multi)
\end{equation}
with notation 
$(h\cdot z)^b=(h\cdot z^1)^{b_1}\cdots (h\cdot z^n)^{b_n}$. 
Then $\xi[h]$ preserves the subspace $V^* \subset \mathbb{C}[V]$. 

\begin{prop}\label{200310.4}
Let $(g,\zeta,q)$ be an admissible triplet. 
\begin{enumerate}
\item For any $h \in G$, the triplet $(hgh^{-1},\zeta,h\cdot q)$ is an admissible triplet. 
$\xi[h]$ induces the isomorphism 
\begin{equation}
\xi[h]:V(g,\zeta,q)(j) \to V(hgh^{-1},\zeta,hq)(j).
\end{equation}
We have 
\begin{eqnarray}
\varphi[hgh^{-1},\zeta,h\cdot q]\circ \xi[h]&=&\xi[h] \circ \varphi[g,\zeta,q],\\
\psi[hgh^{-1},\zeta,h\cdot q] \circ \xi[h]&=&\xi[h] \circ \psi[g,\zeta,q].
\end{eqnarray}
\item 
For any $t \in \mathbb{C}^*$, 
the triplet $(g,\zeta,t\cdot q)$ is an admissible triplet.
For $t \cdot id_V \in GL(V)$ with the identity morphism 
$id_V:V \to V$, 
$\xi[t \cdot id_V]$ preserves the subspace $V^* \subset \mathbb{C}[V]$ and 
$\xi[t \cdot id_V]$ induces the isomorphism 
\begin{equation}
\xi[t \cdot id_V]:V(g,\zeta,q)(j) \to V(g,\zeta,t \cdot q)(j).
\end{equation}
We have 
\begin{eqnarray}
\varphi[g,\zeta,t\cdot q]\circ \xi[t \cdot id_V]&=&\xi[t \cdot id_V] \circ \varphi[g,\zeta,q],\\
\psi[g,\zeta,t\cdot q] \circ \xi[t \cdot id_V]&=&\xi[t \cdot id_V] \circ \psi[g,\zeta,q].
\end{eqnarray}
\item If an integer $r$ satisfies $\mathrm{gcd}(r,d_n)=1$, 
then the triplet $(g^r,\zeta^r,q)$ is also admissible and we have 
\begin{eqnarray}
V(g,\zeta,q)(j)&=&V(g^r,\zeta^r,q)(j),\label{200316.3}\\
\varphi[g,\zeta,q]&=&\varphi[g^r,\zeta^r,q],\label{200316.4}\\
\psi[g,\zeta,q]&=&\psi[g^r,\zeta^r,q]. \label{200316.5}
\end{eqnarray}
\end{enumerate}
\end{prop}
\begin{pf}
We prove (i) and (ii). 
For any $a \in \multi$, we prove $\varphi[hgh^{-1},\zeta,h\cdot q]\circ \xi[h](x^a)
=\xi[h] \circ \varphi[g,\zeta,q](x^a)$. 
\begin{eqnarray*}
&&\varphi[hgh^{-1},\zeta,h\cdot q]\circ \xi[h](x^a)\\
&=&
\varphi[hgh^{-1},\zeta,h\cdot q]((h\cdot x)^a)\\
&=&
\varphi_3[h \cdot q]\circ \varphi_2\circ \varphi_1[h\cdot q]((h\cdot x)^a)\\
&=&
\varphi_3[h \cdot q]\circ \varphi_2(\left(h\cdot x-(h\cdot x)(h\cdot q)\right)^a)\\
&=&
\varphi_3[h \cdot q]\left[
\sum_{b \in \multi}\frac{1}{b!}
\frac{\partial^b \left(h\cdot x-(h\cdot x)(h\cdot q)\right)^a}
{\partial (h\cdot z)^b}(h\cdot q)
\left(
h\cdot z-(h\cdot z)(h\cdot q)
\right)^b
\right]\\
&=&
\varphi_3[h \cdot q]\left[
\sum_{b \in \multi}\frac{1}{b!}
\frac{\partial^b (x-x(q))^a}{\partial z^b}(q)
\left(
h\cdot z-(h\cdot z)(h\cdot q)
\right)^b
\right]\\
&=&
\sum_{b \in \multi}\frac{1}{b!}
\frac{\partial^b (x-x(q))^a}{\partial z^b}(q)
(h\cdot z)^b
\\
&=&
\xi[h] \circ \varphi[g,\zeta,q](x^a). 
\end{eqnarray*}
The other parts are proved in a similar manner, so we omit it. 

We prove (iii). 
Since the triplet  $(g,\zeta,q)$ is admissible, 
we have $g \cdot q=\zeta q$. 
Then $g^r \cdot q=\zeta^r q$. 
Thus the triplet $(g^r,\zeta^r,q)$ is also admissible. 

We show (\ref{200316.3}). 
Let $z^1,\cdots,z^n \in V^*$ be a $g$-homogeneous basis of $V^*$. 
Then for $\one$, $g^r \cdot z^{\alpha}=(\zeta^{d_{\alpha}-1})^r z^{\alpha}
=(\zeta^r)^{d_{\alpha}-1}z^{\alpha}$. 
Then $z^1,\cdots,z^n \in V^*$ be a $g^r$-homogeneous basis of $V^*$. 
Thus we have (\ref{200316.3}). 

Since the morphism $\varphi[g,\zeta,q]$ depends only on the choice of $q$ 
by the proof of Proposition \ref{200310.12}, we have (\ref{200316.4}). 
Then the morphism $\psi[g,\zeta,q]$ depends only on the grading 
$V(g,\zeta,q)(j)$. 
By (\ref{200316.3}), 
we have (\ref{200316.5})
\qed\end{pf}


\begin{defn}
For admissible triplets 
$(g,\zeta,q)$, $(g',\zeta',q')$, we define an equivalence relation 
\begin{equation}
(g,\zeta,q)\sim (g',\zeta',q')
\end{equation}
by 
\begin{equation}
\psi[g,\zeta,q]^{-1}(V^*)=\psi[g',\zeta',q']^{-1}(V^*).
\end{equation} 
\end{defn}
By Proposition \ref{200310.4}, we have the following results.
\begin{cor}
For an admissible triplet $(g,\zeta,q)$ and 
$\forall h \in G$, $\forall t \in \mathbb{C}^*$, $\forall r \in \mathbb{Z}$ 
with $\mathrm{gcd}(d_n,r)=1$, we have 
\begin{eqnarray}
(g,\zeta,q)&\sim& (h\cdot g\cdot h^{-1},\zeta,h\cdot q),\label{200310.5}\\
(g,\zeta,q)&\sim& (g,\zeta,t\cdot q),\label{200310.6}\\
(g^r,\zeta^r,q)&\sim& (g,\zeta,q).\label{200310.7}
\end{eqnarray}
\end{cor}

%

\vskip1cm

\section{Taylor coefficients of the good basic invariants}
Let $G$ be a finite complex reflection group and 
we fix an admissible triplet $(g,\zeta,q)$. 

\begin{defn}\label{200316.10}
A set of basic invariants $x^1,\cdots,x^n$ is 
compatible with a $g$-homogeneous basis $z^1,\cdots,z^n$ of $V^*$ if 
the Jacobian matrix is the identity matrix, i.e. 
$$
\left(
\frac{\p x^{\alpha}}{\p z^{\beta}}(q)
\right)_{1 \leq \alpha,{\beta} \leq n}
=\left(
\delta^{\alpha}_{\beta}
\right)_{1 \leq \alpha,{\beta} \leq n},
$$
where $\delta^{\alpha}_{\beta}$ is the Kronecker's delta. 
\end{defn}

%
%
%
%
\begin{prop}
For a $g$-homogeneous basis $z^1,\cdots,z^n$ of $V^*$, we have the 
following results.
\begin{enumerate}
\item If we put 
\begin{equation}
x^{\alpha}:=\psi[g,\zeta,q]^{-1}(z^{\alpha})\quad
(\one),
\end{equation}
then $x^1,\cdots,x^n$ form a set of basic invariants 
which are good and compatible with a $g$-homogeneous 
basis $z^1,\cdots,z^n$ of $V^*$. 
\item Conversely if $x^1,\cdots,x^n$ are good and 
compatible with a $g$-homogeneous basis $z^1,\cdots,z^n$ of $V^*$, 
then $\psi[g,\zeta,q](x^{\alpha})=z^{\alpha}$ for $\one$. 
\item For any set of basic invariants $x^1,\cdots,x^n$, 
\begin{equation}\label{200310.9}
\frac{\partial^b[x-x(q)]^a}{\partial z^b}(q)=0
\hbox{ if }
d \cdot b \notin \{d\cdot a+d_n j\,|\, j \in \mathbb{Z}_{\geq 0}\} 
\end{equation}
for $a,b \in \multi$. 
\item A set of basic invariants $x^1,\cdots,x^n$ is good if and only if 
\begin{equation}\label{200318.7}
\frac{\p^a x^{\alpha}}{\p z^a}(q)=0\ 
(d_{\alpha}=d\cdot a,\ |a|\geq 2,\ \one). 
\end{equation}
%
%
%
%
\item If a set of basic invariants $x^1,\cdots,x^n$ is good and compatible 
with a $g$-homogeneous basis $z^1,\cdots,z^n$ of $V^*$, then 
for $a,b \in \multi$ satisfying $d \cdot a=d\cdot b$, we have 
\begin{equation}\label{200310.10}
\frac{1}{b!}\frac{\partial^b[x-x(q)]^a}{\partial z^b}(q)=\delta_{a,b}.
\end{equation}
\end{enumerate}
\end{prop}
\begin{pf} As for (i), (ii), they are direct consequences of 
Definition \ref{200316.9}
and 
Definition \ref{200316.10}. 
As for (iii), it is proved in Proposition \ref{200310.12} (ii). 

(iv) By (\ref{200310.14}), we have 
$$
\psi[g,\zeta,q](x^{\alpha})=
\sum_{b \in \multi,\ 
d\cdot b=d_{\alpha}}\frac{1}{b!}
\frac{\partial^b x^{\alpha}}{\partial z^b}(q)z^b.
$$
By the goodness assumption, this must be an element of $V^*$. 
Then the coefficients with $|b| \geq 2$ must be $0$. 

(v) We have $\psi[g,\zeta,q](x^{\alpha})=z^{\alpha}$ for $\one$ 
by (ii). 
Then for any $a \in \multi$, 
\begin{equation}
\psi[g,\zeta,q](x^a)
=\prod_{\gamma=1}^n \psi[g,\zeta,q](x^{\gamma})^{a_i}
=\prod_{\gamma=1}^n (z^{\gamma})^{a_i}
=z^a, 
\end{equation}
and comparing it with (\ref{200310.14}), we have the result. 
\qed\end{pf}

%

\section{The cases of finite real reflection groups}

From now on we shall assume that $G$ is a finite real reflection group, i.e. 
there exists a $G$-stable $\RR$-subspace $V_{\RR}$ of $V$ 
such that the canonical map $V_{\RR}\otimes_{\RR}\CC \to V$ 
is bijective and $G$ is generated by reflections of 
order $2$. We also assume that the action of $G$ on $V$ is irreducible. 


It is known that there exists a set of basic invariants $x^1,\cdots,x^n$ 
and if we denote their degrees by $d_1 \leq \cdots \leq d_n$, 
then they have the following properties (cf. Bourbaki \cite{Bourbaki}).
\begin{enumerate}
\item The action of a Coxeter transformation $g \in G$ on $V$ 
has the eigenvalues 
\begin{equation}
\exp\left(\frac{2\pi\sqrt{-1}(d_1-1)}{d_n}\right),
\cdots,
\exp\left(\frac{2\pi\sqrt{-1}(d_n-1)}{d_n}\right). 
\end{equation}
\item They have a duality
\begin{equation}
d_{\alpha}+d_{n+1-\alpha}=d_1+d_n\quad (\one). \label{200314.2}
\end{equation}
\item We have 
\begin{eqnarray}
d_{n-1}<d_n,\label{200314.1}\\
d_1=2.\label{200314.3}
\end{eqnarray}
\item An eigenvector $v$ of $g$ with the eigenvalue 
$\displaystyle{
\exp\left(\frac{2\pi\sqrt{-1}}{d_n}\right)
}$ is not contained in the reflecting hyperplanes. 
In particular the value of the determinant of the Jacobian 
\begin{equation}\label{200316.11}
\left(
\frac{\p x^{\alpha}}{\p z^{\beta}}(v)
\right)_{1 \leq \alpha,\beta \leq n}
\end{equation}
is nonzero. 
\item There exists a $G$-invariant positive definite symmetric 
$\RR$-bilinear form 
\begin{equation}\label{200316.12}
I_{\RR}:V_{\RR}\times V_{\RR} \to \RR.
\end{equation} 
\end{enumerate}

\subsection{Existence and Uniqueness of a set of good basic invariants}

\begin{prop}\label{230403.2}
For the finite real reflection group $G$, 
there exists an admissible triplet. 
\end{prop}
\begin{pf}
Let $g \in G$ be the Coxeter transformation, 
$v$ be an eigenvector with the eigenvalue 
$\exp\left(\frac{2\pi \sqrt{-1}}{d_n}\right)$. 
Then 
$$
(g,\exp\left(\frac{2\pi \sqrt{-1}}{d_n}\right),v)
$$
satisfies the conditions of the admissibility (i) and (ii) 
by (\ref{200316.11}). 
\qed
\end{pf}
We show the uniqueness of the good invariants. 
\begin{prop} 
For a finite real reflection group $G$, 
any admissible triplets $(g,\zeta,q)$ and $(g',\zeta',q')$ 
satisfy $(g',\zeta',q')\sim (g,\zeta,q)$. 
\end{prop}
\begin{pf}
Take an integer $1 \leq r \leq d_n$ such that $\zeta'=\zeta^r$. 
Since $\zeta,\zeta'$ are primitive $d_n$-th roots of unity, $(d_n,r)=1$. 
Then the triplet $(g^r,\zeta^r,q)=(g^r,\zeta',q)$ is admissible 
by Proposition \ref{200310.4}\hskip0.2mm(iii). 

Compare the admissible triplets $(g^r,\zeta',q)$ and $(g',\zeta',q')$. 
Then we have $g^r=h g' h^{-1}$ for some $h \in G$ 
by a known result (cf. \cite[4.2 Theorem (iv)]{Springer}).  				
By (\ref{200310.5}), we have 
$(g',\zeta',q')\sim (h g' h^{-1},\zeta',h\cdot q')=(g^r,\zeta^r,h\cdot q')$. 

Compare the admissible triplets $(g^r,\zeta^r,h\cdot q')$ and $(g^r,\zeta^r,q)$. 
We see that the dimension of the eigenspace of $g^r$ is $1$ 
because it is the multiplicity 
(
$=\#\{
\beta \in \{1,\cdots,n\}
\,|\,d_{\beta}=d_{\alpha}\,\}$
) 
for $d_{\alpha}$ ($\one$) with $(d_{\alpha},d_n)=1$ and equals $1$ 
for the case $G$ is a finite real reflection group. 
Then using (\ref{200310.6}), we have $(g^r,\zeta^r,h\cdot q') \sim (g^r,\zeta^r,q)$. 

By (\ref{200310.7}), $(g^r,\zeta^r,q)\sim (g,\zeta,q)$. 
Then we have the result. 
\qed\end{pf}
%
%

\begin{cor}\label{200318.1}
For a finite real reflection group $G$, there exists uniquely 
the space of the $\CC$-span of a set of good basic invariants.
\end{cor}

\begin{rmk} 
For an admissible triplet $(g,\zeta,q)$ 
for the finite real reflection group $G$, 
$g$ is not necessarily a Coxeter transformation. 
We give an example. 
Let $G$ be the finite real reflection group of type $H_3$. 
Then the degrees of basic invariants are 
$$
d_1=2,\ 
d_2=6,\ 
d_3=10.
$$
Let $g_0$ be a Coxeter transformation, 
$\zeta$ be $\exp(\frac{2\pi\sqrt{-1}}{d_3})$ and 
$(g_0,\zeta,v)$ be an admissible triplet which is constructed 
in the proof of Proposition \ref{230403.2}. 
By Proposition \ref{200310.4}\hskip0.2mm(iii), $(g_0^3,\zeta^3,v)$ is also an admissible triplet. 
We show that $g_0^3$ is not a Coxeter transformation. 
We compare the eigenvalues of $g_0$ and the ones of $g_0^3$. 
A set of the eigenvalues of $g_0$:
$$
\{
\zeta^{1-d_1}=\zeta^{-1},\ \zeta^{1-d_2}=
\zeta^{-5},\ \zeta^{1-d_3}=\zeta^{-9}\}
$$
and a set of the eigenvalues of $g_0^3$:
$$
\{\zeta^{-3},\ \zeta^{-15},\ \zeta^{-27}\}
$$
do not coincide. Since a Coxeter transformation is unique up to conjugacy, 
we see that $g_0^3$ is not a Coxeter transformation. 
\end{rmk}
%
%

%
%
%

\subsection{Bilinear form}

We take the $G$-invariant positive definite symmetric bilinear form $I_{\RR}$ 
(\ref{200316.12}) and extend it to the $\CC$-bilinear form 
\begin{equation}
I:V \times V \to \CC. 
\end{equation}
Since it is nondegenerate, it induces the $\CC$-bilinear form 
\begin{equation}
I^*:V^* \times V^* \to \CC. 
\end{equation}
This gives 
\begin{equation}\label{240317.1}
I^*:\Omega_{\CC[V]} \otimes_{\CC[V]}\Omega_{\CC[V]} \to \CC[V],
\end{equation}
where $\Omega_{\CC[V]}$ is the module of K\"ahler differentials of $\CC[V]$ over $\CC$. 
It descends to the $\CC[V]^G$-symmetric bilinear form 
\begin{equation}
I^*_G:\Omega_{\CC[V]^G} \otimes_{\CC[V]^G}\Omega_{\CC[V]^G} \to \CC[V]^G
\end{equation}
because $I$ is $G$-invariant, 
where $\Omega_{\CC[V]^G}$ is the module of K\"ahler differentials of $\CC[V]^G$ over $\CC$. 

For a set of basic invariants $x^1,\cdots,x^n$, we have 
\begin{equation}
I^*_G(dx^{\alpha},dx^{\beta})=
\sum_{\gamma_1,\gamma_2=1}^{n}
\frac{\p x^{\alpha}}{\p z^{\gamma_1}}
\frac{\p x^{\beta}}{\p z^{\gamma_2}}
I^*(z^{\gamma_1},z^{\gamma_2}) 
\end{equation}
for $1 \leq \alpha,\beta \leq n$. 

\subsection{Good basic invariants and Bilinear form}
From now on we fix an admissible triplet $(g,\zeta,q)$ 
for a finite real reflection group $G$ which acts on $V$. 

For a set of basic invariants $x^1,\cdots,x^n$, we have 
\begin{equation}\label{200310.8}
x^n(q) \neq 0
\end{equation}
because if $x^n(q)=0$, 
then $x^{\alpha}(q)=0$ $(1 \leq \alpha \leq n)$ 
by (\ref{200310.2}) and (\ref{200314.1}), which contradicts $q \neq 0$. 
\begin{thm}\label{230403.3}
Let $(g,\zeta,q)$ be the admissible triplet. 
We assume that a $g$-homogeneous basis $z^1,\cdots,z^n$ of $V^*$ satisfies 
\begin{equation}\label{200310.11}
I^*(z^{\alpha},z^{\beta})=\delta_{\alpha+\beta,n+1} \quad
(1 \leq \alpha,\beta \leq n). 
\end{equation}
Let $x^1,\cdots,x^n$ be a set of good basic invariants 
compatible with this basis $z^1,\cdots,z^n$ of $V^*$. 
Then by Taylor coefficients  
\begin{equation}\label{200318.5}
\frac{\partial^a x^{\alpha}}
{\partial z^a}(q)\quad
(\one,\,a \in \multi,\, d\cdot a=d_{\alpha}+d_n), 
\end{equation}
any 
$I^*_G(dx^{\alpha},dx^{\beta})$ $(\alpha,\beta=1,\cdots,n)$ 
is written as follows:
\begin{eqnarray}\label{240317.2}
&&I^*_G(dx^{\alpha},dx^{\beta})\label{200318.6}\\
&=&\delta_{\alpha+\beta,n+1}\frac{x^{n}}{x^{n}(q)}
+\sum_{\substack{b=(b_1,\cdots,b_n),\\ b_n=0,\\ d\cdot b=d_{\alpha}+d_{\beta}-2}}
\frac{1}{b!}
\left[
\frac{\p^b}{\p z^b}
\left(
\frac{\partial x^{\alpha}}{\partial z^{\beta*}}
+
\frac{\partial x^{\beta}}{\partial z^{\alpha*}}
\right)
\right](q)\,
x^b,\nonumber
\end{eqnarray}
where $\alpha*=n+1-\alpha \quad (\one)$. 
\end{thm}
\begin{pf}
For any $\alpha,\beta\ (\two)$, $I^*_G(dx^{\alpha},dx^{\beta})$ 
is represented as  
\begin{equation}\label{200309.3}
I^*_G(dx^{\alpha},dx^{\beta})
=\sum_{\substack{a=(a_1,\cdots,a_n) \in \multi,\\
 a_n=0,\\ d\cdot a=d_{\alpha}+d_{\beta}-d_n-2}}
A_a^{\alpha,\beta}
x^ax^{n}
+
\sum_{\substack{b=(b_1,\cdots,b_n)\in \multi,\\
 b_n=0,\\ d\cdot b=d_{\alpha}+d_{\beta}-2}}
B_b^{\alpha,\beta}
x^b
\end{equation}
for $A_a^{\alpha,\beta},B_b^{\alpha,\beta} \in \mathbb{C}$ 
by the degree of $\CC[V]^G$. 

By taking higher order derivatives of the both sides of (\ref{200309.3}) 
with respect to $z^1,\cdots,z^n$ and evaluating them at $q$, 
we determine $A_a^{\alpha,\beta},B_b^{\alpha,\beta}$ 
in the following lemmas. 
\begin{lem}
For the cases $d_{\alpha}+d_{\beta} \leq d_n+2$, 
\begin{equation}
A_a^{\alpha,\beta}=0\ \mbox{if $d_{\alpha}+d_{\beta} <d_n+2$ or $a \neq 0$}.
\end{equation}
\end{lem}
\begin{pf}
The multi-indices $a \in \multi$ satisfying 
$d\cdot a=d_{\alpha}+d_{\beta}-d_n-2\leq 0$ 
are only $a=(0,\cdots,0) \in \multi$. Thus we have 
the results. 
\qed
\end{pf}
\begin{lem}
We have 
\begin{equation}
(\mathrm{RHS} \hbox{ of }(\ref{200309.3}))(q)=
\begin{cases}
A_0^{\alpha,\beta}x^{n}(q)\ &\mbox{if $d_{\alpha}+d_{\beta} = d_n+2$},\\
0\ &\mbox{if $d_{\alpha}+d_{\beta} \neq d_n+2$}.
\end{cases}
\end{equation}
\end{lem}
\begin{pf}
By 
$(A_a^{\alpha,\beta}x^ax^{n})(q)=0$ if $a \neq 0$ and 
$(B^{\alpha,\beta}_bx^b)(q)=0$ 
by $x^{\alpha}(q)=0\ (1 \leq \alpha \leq n-1)$ 
which are shown in (\ref{200310.2}), we have the results. 
\qed
\end{pf}
\begin{lem}
We have 
\begin{equation}
(\mathrm{LHS} \hbox{ of }(\ref{200309.3}))(q)
=\delta_{\alpha+\beta,n+1}. 
\end{equation}
\end{lem}
\begin{pf}
We have 
\begin{eqnarray*}
\hbox{(LHS of (\ref{200309.3}))}(q)
&=&\sum_{\gamma_1,\gamma_2=1}^n 
\frac{\partial x^{\alpha}}{\partial z^{\gamma_1}}(q)
\frac{\partial x^{\beta}}{\partial z^{\gamma_2}}(q)
I^*(z^{\gamma_1},z^{\gamma_2})\nonumber \\
&=&\sum_{\gamma_1,\gamma_2=1}^n 
\delta_{\alpha,\gamma_1}
\delta_{\beta,\gamma_2}
\delta_{\gamma_1+\gamma_2,n+1}\\
&=&\delta_{\alpha+\beta,n+1}. 
\end{eqnarray*}
\qed
\end{pf}
\begin{lem}
For the cases $d_{\alpha}+d_{\beta}>d_n+2$, 
take any multi-index $c=(c_1,\cdots,c_n)\in \multi$ such that 
$c_n=0,\ d\cdot c=d_{\alpha}+d_{\beta}-d_n-2$, we have 
\begin{equation}\label{200309.4}
\left[
\frac{1}{c!}\frac{\partial^c}{\partial z^c}
(\mathrm{RHS} \hbox{ of }(\ref{200309.3}))
\right]
(q)=A_c^{\alpha,\beta}x^{n}(q). 
\end{equation}
\end{lem}
\begin{pf}
For $A^{\alpha,\beta}_a x^a x^{n}$ with any 
$a=(a_1,\cdots,a_n)\in \multi,\ a_n=0,\ d\cdot a=d_{\alpha}+d_{\beta}-d_n-2$, 
\begin{equation}\label{200313.1}
\frac{1}{c!}\frac{\partial^c(A^{\alpha,\beta}_a x^ax^{n})}{\partial z^c}(q)
\end{equation}
is a linear combination of 
\begin{equation}\label{200313.4}
A^{\alpha,\beta}_a
\left[
\frac{\partial^{c'}(x^a)}{\partial z^{c'}}(q)
\right]
\left[
\frac{\partial^{c''}(x^{n})}{\partial z^{c''}}(q)
\right]\quad
(c',\,c'' \in \multi)
\end{equation}
with $c'+c''=c$. 
Since $c''_n=0$ and 
$d\cdot c'' \leq d\cdot c =d_{\alpha}+d_{\beta}-d_n-2<d_n$, 
(\ref{200313.4}) must be $0$ by (\ref{200310.9}) if $c'' \neq 0$. 
Then (\ref{200313.1}) equals 
\begin{equation*}
A^{\alpha,\beta}_a \frac{1}{c!}\frac{\partial^c(x^a)}{\partial z^c}(q)x^{n}(q). 
\end{equation*}
Since $d\cdot a=d\cdot c$ and $x^{a}=(x-x(q))^a$ by $a_n=0$, we have 
\begin{equation*}
A^{\alpha,\beta}_a 
\frac{1}{c!}\frac{\partial^c(x^a)}{\partial z^c}(q)x^{n}(q)
=A^{\alpha,\beta}_a \delta_{a,c} x^{n}(q)
\end{equation*}
by (\ref{200310.10}).

For $B^{\alpha,\beta}_b x^b$ with any 
$b=(b_1,\cdots,b_n)\in \multi,\ b_n=0,\ d\cdot b=d_{\alpha}+d_{\beta}-2$, 
we have 
\begin{equation*}
\frac{1}{c!}\frac{\partial^c (B^{\alpha,\beta}_b x^b)}{\partial z^c}(q)=0
\end{equation*}
by $d\cdot c<d\cdot b$ and (\ref{200310.9}). 

Then we have the equation (\ref{200309.4}). 
\qed
\end{pf}

\begin{lem}
For the cases $d_{\alpha}+d_{\beta}>d_n+2$, 
take any multi-index $c=(c_1,\cdots,c_n)\in \multi$ such that 
$c_n=0,\ d\cdot c=d_{\alpha}+d_{\beta}-d_n-2$, we have 
\begin{equation}\label{200309.5}
\left[
\frac{1}{c!}\frac{\partial^c}{\partial z^c}
(\mathrm{LHS} \hbox{ of }(\ref{200309.3}))
\right]
(q)=0. 
\end{equation}
\end{lem}
\begin{pf}
The LHS of (\ref{200309.5}) is a linear combination of the products:
\begin{equation}\label{210721.2}
\left[
\frac{\partial^{c'}}{\partial z^{c'}}
\frac{\partial x^{\alpha}}{\partial z^{\gamma_1}}
\right](q)
\left[
\frac{\partial^{c''}}{\partial z^{c''}}
\frac{\partial x^{\beta}}{\partial z^{\gamma_2}}
\right](q)
I^*(z^{\gamma_1},z^{\gamma_2})\quad
(c',\,c'' \in \multi)
\end{equation}
with $c'+c''=c$. 
We assume that (\ref{210721.2}) is nonzero for some $c',c''$. 
Then there exists non-negative integers $n_1,n_2$ such that 
\begin{equation*}
d \cdot c'+d_{\gamma_1}=d_{\alpha}+n_1d_n, \quad
d \cdot c''+d_{\gamma_2}=d_{\beta}+n_2d_n
\end{equation*}
by (\ref{200310.9}). We also have 
\begin{equation*}
d_{\gamma_1}+d_{\gamma_2}=d_n+2
\end{equation*}
by (\ref{200310.11}). 
Combining these equalities, we have $n_1+n_2=0$. 
Then we have $n_1=n_2=0$. 
Then we have 
\begin{equation*}
\left[
\frac{\partial^{c'}}{\partial z^{c'}}
\frac{\partial x^{\alpha}}{\partial z^{\gamma_1}}
\right](q)\neq 0
\hbox{ for }
d \cdot c'+d_{\gamma_1}=d_{\alpha}.
\end{equation*}
Since $x^{\alpha}$ satisfies the condition (\ref{200318.7}), 
we have $c'=0$ and $\alpha=\gamma_1$. 
Also we have $c''=0$. Then $c=c'+c''=0$ which contradicts $d\cdot c>0$. 
Therefore we have the equation (\ref{200309.5}). 
\qed
\end{pf}
By these Lemmas, for any $\alpha,\beta,c$ with $\two,\ c \in \multi$, we obtain 
\begin{eqnarray}\label{200313.2}
A^{\alpha,\beta}_c=
\begin{cases}
\delta_{\alpha+\beta,n+1}/x^{n}(q) 
&\mbox{if $d_{\alpha}+d_{\beta}=d_n+2$ and $c=0$},\\
0
&\mbox{otherwise},
\end{cases}
\end{eqnarray}
where we used $x^{n}(q)\neq 0$ by (\ref{200310.8}). 
\begin{lem}
For any multi-index $c \in \multi$ with 
$c_n=0,\ d\cdot c=d_{\alpha}+d_{\beta}-2$, we have 
\begin{equation}\label{200309.6}
\left[
\frac{1}{c!}\frac{\partial^{c}}{\partial z^{c}}
(\mathrm{RHS} \hbox{ of }(\ref{200309.3}))
\right]
(q)=B_{c}^{\alpha,\beta}. 
\end{equation}
\end{lem}
\begin{pf}
For $A^{\alpha,\beta}_a x^a x^{n}$ with any 
$a=(a_1,\cdots,a_n) \in \multi,\ a_n=0,\ d\cdot a=d_{\alpha}+d_{\beta}-d_n-2$, 
we show 
\begin{equation}\label{200318.9}
\frac{1}{c!}\frac{\partial^c(A^{\alpha,\beta}_a x^ax^{n})}{\partial z^c}(q).
\end{equation}
We have only to prove it for the cases $A^{\alpha,\beta}_a \neq 0$. 
Then $d_{\alpha}+d_{\beta}=d_n+2$ and $a=0$ by 
(\ref{200313.2}). 
Then by $d\cdot c=d_n$ and $c_n=0$, we have (\ref{200318.9}). 

For $B^{\alpha,\beta}_b x^b$ with any 
$b=(b_1,\cdots,b_n)\in \multi,\ b_n=0,\ d\cdot b=d_{\alpha}+d_{\beta}-2$, 
we have 
\begin{equation*}
\frac{1}{c!}\frac{\p^{c} (B^{\alpha,\beta}_{b} x^{b})}{\p z^{c}}(q)
=B^{\alpha,\beta}_{b}\delta_{c,b}
=B^{\alpha,\beta}_{c}
\end{equation*}
by $d \cdot c=d\cdot b$ and (\ref{200310.10}). 
Thus we have (\ref{200309.6}). 
\qed
\end{pf}
\begin{lem}
For any multi-index $c \in \multi$ with 
$c_n=0,\ d\cdot c=d_{\alpha}+d_{\beta}-2$, we have 
\begin{equation}\label{200309.7}
\left[
\frac{1}{c!}\frac{\partial^c}{\partial z^c}
(\mathrm{LHS} \hbox{ of }(\ref{200309.3}))
\right]
(q)=\frac{1}{c!}
\left[
\left(
\frac{\partial^{c}}{\partial z^{c}}
\frac{\partial x^{\alpha}}{\partial z^{n+1-\beta}}
\right)(q)
+
\left(
\frac{\partial^{c}}{\partial z^{c}}
\frac{\partial x^{\beta}}{\partial z^{n+1-\alpha}}
\right)
\right](q). 
\end{equation}
\end{lem}
\begin{pf}
For the LHS of (\ref{200309.7}), it is a linear combination of the products:
\begin{equation}\label{210721.3}
\left[
\frac{\partial^{c'}}{\partial z^{c'}}
\frac{\partial x^{\alpha}}{\partial z^{\gamma_1}}
\right](q)
\left[
\frac{\partial^{c''}}{\partial z^{c''}}
\frac{\partial x^{\beta}}{\partial z^{\gamma_2}}
\right](q)
I^*(z^{\gamma_1},z^{\gamma_2})\quad
(c',\,c'' \in \multi)
\end{equation}
with $c'+c''=c$. 
We assume that (\ref{210721.3}) is nonzero for some $c',c''$. 
Then there exists non-negative integers $n_1,n_2$ such that 
$d \cdot c'+d_{\gamma_1}=d_{\alpha}+n_1d_n$, 
$d \cdot c''+d_{\gamma_2}=d_{\beta}+n_2d_n$ by (\ref{200310.9}) and 
$d_{\gamma_1}+d_{\gamma_2}=d_n+2$ by (\ref{200310.11}). 
Combining these equalities, we have $n_1+n_2=1$. 
Then we have $(n_1,n_2)=(1,0),(0,1)$. 
For the case $(n_1,n_2)=(1,0)$, 
$d\cdot c''+d_{\gamma_2}=d_{\beta}$. 
Since $x^{\beta}$ satisfies (\ref{200310.10}), 
$c''$ must be $0$ and $\gamma_2=\beta$. 
Then we have $c'=c$ and $\gamma_1=n+1-\beta$. 
For the case $(n_1,n_2)=(0,1)$, 
$c'$ must be $0$ and $\gamma_1=\alpha$. 
Then we have $c''=c$ and $\gamma_2=n+1-\alpha$. 

Then we have 
\begin{eqnarray*}
&&c! \times [\hbox{LHS of (\ref{200309.7})}]\nonumber \\
&&=
\left[
\frac{\partial^{c}}{\partial z^{c}}
\frac{\partial x^{\alpha}}{\partial z^{n+1-\beta}}
\right](q)
\frac{\partial x^{\beta}}{\partial z^{\beta}}
(q)
I^*(z^{n+1-\beta},z^{\beta})
+
\frac{\partial x^{\alpha}}{\partial z^{\alpha}}
(q)
\left[
\frac{\partial^{c}}{\partial z^{c}}
\frac{\partial x^{\beta}}{\partial z^{n+1-\alpha}}
\right](q)
I^*(z^{\alpha},z^{n+1-\alpha})
\nonumber
\\
&&=
\left[
\frac{\partial^{c}}{\partial z^{c}}
\frac{\partial x^{\alpha}}{\partial z^{n+1-\beta}}
\right](q)
+
\left[
\frac{\partial^{c}}{\partial z^{c}}
\frac{\partial x^{\beta}}{\partial z^{n+1-\alpha}}
\right](q).
\end{eqnarray*}
Then we have the equation (\ref{200309.7}). 
\qed
\end{pf}
By (\ref{200309.6}) and (\ref{200309.7}), we have 
\begin{equation}
B_c^{\alpha,\beta} 
=
\frac{1}{c!}
\left[
\frac{\p^c}{\p z^c}
\left(
\frac{\partial x^{\alpha}}{\partial z^{n+1-\beta}}
+
\frac{\partial x^{\beta}}{\partial z^{n+1-\alpha}}
\right)
\right](q)\,. 
\end{equation}

\qed\end{pf}
\section{Good basic invariants and Frobenius structure}
\subsection{Euler field}

We introduce the Euler field, 
the space of lowest degree part of the derivatives 
and the $G$-invariant bilinear form. 

Let $x^1,\cdots,x^n$ be a set of basic invariants with 
degrees $d_1 \leq \cdots \leq d_{n-1}<d_n$ (see (\ref{200314.1})). 

We define the Euler field $E$ by 
\begin{eqnarray}
E&:=&\sum_{\alpha=1}^n \frac{d_{\alpha}}{d_n}x^{\alpha}\frac{\p}{\p x^{\alpha}}
: \Omega_{\CC[V]^G} \to \CC[V]^G,
\end{eqnarray}
which does not depend on the choice of a set of basic invariants. 

Let $\Der$ be the module of $\CC$-derivations of $\CC[V]^G$. 
It has the grading  
\begin{equation}
\Der=\bigoplus_{j \in \ZZ}\Der (j),\quad
\frac{\p}{\p x^{\alpha}} \in \Der (-d_{\alpha})\quad
(\one)
\end{equation}
induced by the grading of $\CC[V]^G=\bigoplus_{j \in \ZZ}S(j)$ 
and we see that the dimension of the lowest degree part is 
\begin{equation}
\mathrm{dim}_{\CC} \Der (-d_{n})=1
\end{equation}
by (\ref{200314.1}). 

\subsection{Frobenius structure}
The Frobenius structure on $\CC[V]^G$ is constructed by Saito 
\cite{finite} and Dubrovin \cite{Dubrovin} 
(see also \cite{Hertling}). 
\begin{thm}\label{200318.2}
$($
Saito \cite{finite}, Dubrovin \cite{Dubrovin}
$)$ 
\begin{enumerate}
\item 
There exist a ${\CC[V]^G}$-nondegenerate symmetric 
bilinear form $($called the metric$)$ 
$J:\Der \otimes_{\CC[V]^G} \Der \to {\CC[V]^G}$, 
a ${\CC[V]^G}$-symmetric bilinear form $($called the multiplication$)$ 
$\circ:\Der \otimes_{{\CC[V]^G}}\Der \to \Der$ and 
a field $e \in Der(\CC[V]^G)$ 
subject to the following conditions:
\begin{enumerate}
\item the metric is invariant under the multiplication, 
i.e. $J(X\circ Y,Z)=J(X,Y \circ Z)$ 
for any vector fields $X,Y,Z:\Omega_{\CC[V]^G} \to {\CC[V]^G}$, 
\item $($potentiality$)$ the 
$(3,1)$-tensor $\nabla \circ$ is symmetric
(where $\nabla$ is the Levi-Civita connection of the metric), 
i.e. 
$\nabla_X(Y \circ Z)
-Y\circ \nabla_X (Z)
-\nabla_Y(X \circ Z)
+X\circ \nabla_Y (Z)
-[X,Y]\circ Z=0$, 
for any vector fields $X,Y,Z:\Omega_{\CC[V]^G} \to {\CC[V]^G}$, 
\item the metric $J$ is flat, 
\item $e$ is a unit field for $\circ$ and it is flat, i.e. $\nabla e=0$,
\item the Euler field $E$ satisfies 
$Lie_E(\circ)=1\cdot \circ$, and 
$Lie_E(J)=(2-\frac{d_n-2}{d_n})\cdot J$, 
\item the intersection form coincides with the bilinear form $I^*_G$:
$J(E,J^*(\omega)\circ J^*(\omega'))=I^*_G(\omega,\omega')$ 
for 1-forms $\omega,\omega' \in \Omega_{\CC[V]^G}$ and 
$J^*:\Omega_{\CC[V]^G} \to \Der$ is the isomorphism induced by 
the dual metric $J^*$ of $J$. 
\end{enumerate}
\item Let $(J,\circ,e)$ be a Frobenius structure 
satisfying the conditions in $(\mathrm{i})$. 
Then $e \in Der(\CC[V]^G)(-d_{n}) \setminus\{0\}$. 
Conversely for any element $\widetilde{e} \in 
Der(\CC[V]^G)(-d_{n}) \setminus \{0\}$, 
there exists uniquely a Frobenius structure 
$(\widetilde{J},\widetilde{\circ},\widetilde{e})$ satisfying 
the conditions in $(\mathrm{i})$. 
The Frobenius structure 
$(\widetilde{J},\widetilde{\circ},\widetilde{e})$ 
is written as 
$(\widetilde{J},\widetilde{\circ},\widetilde{e})
=(c^{-1}J,c^{-1}\circ,ce)$ for some $c \in \CC^{\times}$. 
\end{enumerate}
\end{thm}
The metric $J$ could be constructed from $I^*_G$ and 
$e$ as follows. 
\begin{prop}\label{200323.1}
For 1-forms $\omega,\omega' \in \Omega_{\CC[V]^G}$, 
we have 
\begin{equation}
J^*(\omega,\omega')=(Lie_{e}(I^*_G))(\omega,\omega')
\end{equation}
for the dual metric $J^*$ of $J$. 
\end{prop}
\begin{pf}
By combining the results 
$Lie_e(J)=0$, $Lie_e(\circ)=0$ and $Lie_e(E)=e$ (cf.\cite[p146]{Hertling})   
with the Lie derivative of the both sides of the equation 
$$
J(E,J^*(\omega)\circ J^*(\omega'))=I^*_G(\omega,\omega')
$$
in Theorem \ref{200318.2}\hskip0.2mm(i)(f) with respect to the unit $e$, 
we have the result.
\qed
\end{pf}
Let $\nabla$ be a connection introduced in Theorem \ref{200318.2}. 
By Theorem \ref{200318.2}\hskip0.2mm(ii), 
the metric $J$ of the Frobenius structure satisfying 
conditions in Theorem \ref{200318.2}\hskip0.2mm(i) 
is unique up to a constant factor. 
Then $\nabla$ and the notion of {\it flatness} 
do not depend on the choice of the Frobenius structures 
in Theorem \ref{200318.2}. 
\begin{defn}
A set of basic invariants $x^1,\cdots,x^n$ is called flat with respect to 
the Frobenius structure if 
\begin{equation}
\nabla dx^{\alpha}=0\quad
(\one).
\end{equation}
\end{defn}
\vskip1cm

\subsection{Frobenius structure via flat basic invariants }
We give a description of the multiplication and the metric 
with respect to the set of flat basic invariants. 
\begin{prop}\label{200318.4}
A set of basic invariants $x^1,\cdots,x^n$ with degrees 
$d_1 \leq \cdots \leq d_{n-1}<d_n$ is flat with respect to the 
Frobenius strucuture $(J,\circ,e)$ in Theorem $\ref{200318.2}$ if and only if 
\begin{equation}\label{200318.3}
\eta^{\alpha,\beta}:=eI^*_G(dx^{\alpha},dx^{\beta})\quad
(1 \leq \alpha,\beta \leq n)
\end{equation}
are all elements of $\CC$. 
If a set of basic invariants $x^1,\cdots,x^n$ is flat, 
then the metric $J$ is described by 
\begin{equation}\label{200323.2}
\left(
\eta_{\alpha,\beta}
\right)_{1 \leq \alpha,\beta \leq n}
:=
\left(
J
\left(
\p_{\alpha},\p_{\beta}
\right)
\right)
_{1 \leq \alpha,\beta \leq n}
=
\left(
\eta^{\alpha,\beta}
\right)^{-1}_{1 \leq \alpha,\beta \leq n}
\end{equation}
and the structure constants $C_{\alpha,\beta}^{\gamma}$ of the multiplication 
defined by 
\begin{equation}
\p_{\alpha}{\circ} \p_{\beta}
=\sum_{\gamma=1}^{n}
C_{\alpha,\beta}^{\gamma}
\p_{\gamma} \quad
(1 \leq \alpha,\beta \leq n)
\end{equation}
are described by 
\begin{equation}
C_{\alpha,\beta}^{\gamma}=
\sum_{\alpha',\beta'=1}^{n}
\eta_{\alpha,\alpha'}
\eta_{\beta,\beta'}
\p^{\gamma}
\left(
\frac{d_n}{d_{\alpha'}+d_{\beta'}-2}I^*_G(dx^{\alpha'},dx^{\beta'})
\right) 
\end{equation}
for $1 \leq \alpha,\beta \leq n$, 
where we denote 
\begin{eqnarray}
\p_{\alpha}=\frac{\p}{\p x^{\alpha}},\quad
\p^{\alpha}=\sum_{\alpha'=1}^{n}
\eta^{\alpha,\alpha'}\frac{\p}{\p x^{\alpha'}}\quad (\one).
\end{eqnarray}
\end{prop}

%
%

\begin{pf}
By Proposition \ref{200323.1}, the dual metric of the metric 
of the Frobenius strucrure is constructed 
from the unit $e$ and $I^*_G$ by (\ref{200318.3}). 

For the construction of the multiplication from $I^*_G$, 
we remind the notion of the Frobenius potential (see Dubrovin \cite{Dubrovin}). 

The Frobenius potential $F$ is defined by the relation 
\begin{equation}
C_{\alpha,\beta}^{\gamma}=\p_{\alpha}\p_{\beta}\p^{\gamma}F
\quad
(1\leq \alpha,\beta,\gamma \leq n)
\end{equation}
with the structure constants 
$C_{\alpha,\beta}^{\gamma}$ of the product 
and it is related with $I^*_G$ as 
\begin{equation}
I^*_G(dx^{\alpha},dx^{\beta})
=
\frac{d_{\alpha}+d_{\beta}-2}{d_n}
\p^{\alpha}
\p^{\beta}F\quad
(1 \leq \alpha,\beta \leq n). 
\end{equation}
Then for any $\alpha,\beta,\gamma\ (1 \leq \alpha,\beta,\gamma \leq n)$, 
we have 
\begin{eqnarray}
C_{\alpha,\beta}^{\gamma}
&=&
\p_{\alpha}\p_{\beta}\p^{\gamma}F \nonumber\\
&=&
\sum_{\alpha',\beta'=1}^{n}
\eta_{\alpha,\alpha'}
\eta_{\beta,\beta'}
\partial^{\gamma}\p^{\alpha'}\p^{\beta'}F\nonumber\\
&=&
\sum_{\alpha',\beta'=1}^{n}
\eta_{\alpha,\alpha'}
\eta_{\beta,\beta'}
\p^{\gamma}
\left(
\frac{d_n}{d_{\alpha'}+d_{\beta'}-2}
I^*_G(dx^{\alpha'},dx^{\beta'})
\right).
\end{eqnarray}
Then we have the results. 
\qed
\end{pf}
\subsection{Good basic invariants and Frobenius structure}

\begin{cor}
\begin{enumerate}
\item Let $x^1,\cdots,x^n$ be the same as in Theorem $\ref{230403.3}$. 
Let $J$ be a metric and $\circ$ be a multiplication 
of a unique Frobenius structure with the unit 
$e=x^n(q)\frac{\p}{\p x^n}$ in Theorem $\ref{200318.2}$. 
Then the metric $J$ and the structure constants of the multiplication 
$C_{\alpha,\beta}^{\gamma}$ 
$(1 \leq \alpha,\beta,\gamma \leq n)$ are 
\begin{eqnarray}
&&J(\frac{\p}{\p x^{\alpha}},\frac{\p}{\p x^{\beta}})
=\delta_{\alpha+\beta,n+1} \quad 
(1 \leq \alpha, \beta \leq n),\label{200323.3}\\
&&C_{\alpha,\beta}^{\gamma}=
\frac{\p }{\p x^{\gamma*}}
\left(
\frac{d_n}{d_{\alpha*}+d_{\beta*}-2}I^*_G(dx^{\alpha*},dx^{\beta*})
\right), \label{240317.3}
\end{eqnarray}
which are all written 
by Taylor coefficients $(\ref{200318.5})$ 
by $(\ref{200318.6})$. 
\item If a set of basic invariants is good $($which is independent of the 
choice of the admissible triplet by Corollary $\ref{200318.1})$, 
then it is flat with respect to the Frobenius structure of Thoerem $\ref{200318.2}$. 
\item The space $\mathrm{Spec}\,\CC[V]$ has a metric induced by the dual 
metric $I^*$ $(\ref{240317.1})$. The space $\mathrm{Spec}\,\CC[V]^G$ has a metric $J$. 
Then $\psi[g,\zeta,q]:\CC[V]^G \simeq \CC[V]$ 
gives the isometry with respect to these metric structures. 
\end{enumerate}
\end{cor}
\begin{pf}
We prove (i). 
Let $x^1,\cdots,x^n$ be the same in Theorem \ref{230403.3}. 
Since $x^n(q) \neq 0$ by (\ref{200310.8}), 
$e=x^n(q)\frac{\p}{\p x^n} \in \Der(-d_n) \setminus \{0\}$. 
By Theorem \ref{200318.2}\hskip0.2mm(ii), we have a unique Frobenius structure 
with the unit $e=x^n(q)\frac{\p}{\p x^n}$. 
By Theorem \ref{230403.3}, we have 
\begin{equation}
eI^*_G(dx^{\alpha},dx^{\beta})=\delta_{\alpha+\beta,n+1}\quad 
(1 \leq \alpha, \beta \leq n). 
\end{equation}
By Proposition \ref{200318.4}, a set of $x^1,\cdots,x^n$ is flat 
and we have (\ref{200323.3}) and (\ref{240317.3}). 
By (\ref{240317.2}) in Theorem \ref{230403.3}, 
(\ref{240317.3}) are all written 
by Taylor coefficients (\ref{200318.5}).

(ii) is a direct consequence of (i). 
(iii) is a direct consequence of (\ref{200310.11}), 
(\ref{200323.3}) 
and $\psi[g,\zeta,q](x^{\alpha})=z^{\alpha}$ 
for $\one$. 
\qed\end{pf}

\end{document}